\documentclass[12pt,a4paper]{article} 
\usepackage{amssymb}
\usepackage{times} 
\usepackage[all]{xy}
\usepackage{hyperref}
\hypersetup{pdfstartview=FitH}
\hypersetup{pdfpagemode=FitWidth} 
\hypersetup{bookmarksnumbered}
\newcommand{\nc}{\newcommand} 

\nc{\op}[1]{\mathop{\mathsf{#1}}\nolimits}

\nc{\bb}{\bigskip}
\nc{\C}{\mathbb{C}}
\nc{\CE}{\op{CE}}
\nc{\cl}{\centerline} 
\nc{\comp}{\op{comp}}\nc{\ext}{\op{ext}}
\nc{\D}{\mathsf{D}}
\nc{\epi}{\twoheadrightarrow}
\nc{\Ext}{\op{Ext}}
\nc{\g}{\langle}\nc{\h}{\rangle}
\nc{\G}{\Gamma}
\nc{\go}{\mathfrak} 
\nc{\Hom}{\op{Hom}}
\nc{\HH}{\mathsf{H}}
\nc{\Id}{\op{Id}} 
\nc{\incl}{\hookrightarrow}
\nc{\ind}{\hskip 2em\relax}
\nc{\iso}{\stackrel{\sim}{\rightarrow}}
\nc{\aitem}{\begin{itemize}}\nc{\zitem}{\end{itemize}}
\nc{\K}{\mathsf{K}}
\nc{\Ker}{\op{Ker}}
\nc{\mc}{\mathcal}
\nc{\Mod}{\op{Mod}}
\nc{\mono}{\rightarrowtail} 
\nc{\N}{\mathbb{N}}
\nc{\ov}{\overline}
\nc{\pt}{\bullet}
\nc{\QC}{\op{QC}}
\nc{\RHom}{\op{RHom}}
\nc{\Ri}{\mathsf{R}}
\nc{\strincl}{\subsetneqq}
\nc{\Supp}{\op{Supp}}
\nc{\then}{\Longrightarrow}
\nc{\Tot}{\mathsf{Tot}}
\nc{\Wedge}{\mathop{\bigwedge}\nolimits} 
\nc{\Z}{\mathbb{Z}}

\newtheorem{ttt}{Theorem}
\newtheorem{cor}[ttt]{Corollary}

\newtheorem{lem}[ttt]{Lemma}
\newtheorem{ppp}[ttt]{Proposition}

\begin{document}

\cl{\Huge Grothendieck categories}\bb

\cl{\Huge and support conditions}\bb\bb

We give examples of pairs $(\mc{G}_1,\mc{G}_2)$ where $\mc{G}_1$ is 
a Grothendieck category and $\mc{G}_2$ a full Grothendieck subcategory of 
$\mc{G}_1$, the inclusion $\mc{G}_2\incl\mc{G}_1$ being denoted 
$\iota$, for which $\Ri^+\iota:\D^+\mc{G}_2\to\D^+\mc{G}_1$ (or even 
$\Ri\iota:\D\mc{G}_2\to\D\mc{G}_1$) is a a full embedding\footnote{%
The categories $\mc{G}_1$ and $\mc{G}_2$ will come under various names, 
but the inclusion will always be denoted by $\iota$.}. 
This yields generalizations of some results of Bernstein and Lunts, 
and of Cline, Parshall and Scott. To wit, Theorem~\ref{bl} (resp. 
Theorem~\ref{cpso}, resp. Theorem~\ref{cpst} 
and Corollary~\ref{cpsc}) below strengthen Theorem 17.1 in Bernstein and 
Lunts \cite{bl} (resp. Example 3.3.c and Theorem 3.9.a of Cline, Parshall 
and Scott \cite{cps2}, resp. Theorem 3.1 and Proposition 3.6 of Cline, 
Parshall and Scott \cite{cps1}). \bb 

\ind We work in the axiomatic system defined by Bourbaki in 
\cite{bkiens}. We postulate in addition the existence of an 
uncountable universe $\mc{U}$ in the sense of Bourbaki \cite{bkisga}. 
All categories are $\mc{U}$-categories. \bb 

\ind By Alonso Tarr\'{\i}o, Jerem\'{\i}as L\'opez and Souto 
Salorio \cite{ajs}, Theorem 5.4, or by Serp\'e \cite{serpe} Theorem 3.13, 
(or more simply by Spaltenstein \cite{spalt}, proof of Theorem 4.5), the 
functor\footnote%
{An example of category for which $\RHom$ can be explicitly described 
is given in Appendix~1.}
$\RHom_{\mc{G}_i}$ is defined on the whole of 
$\D\mc{G}_i^{op}\times\D\mc{G}_i$. --- Consider the following conditions. 
\bb

(R) : For all $V,W\in\D\mc{G}_2$ the complexes $\RHom_{\mc{G}_1}(V,W$) and 
$\RHom_{\mc{G}_2}(V,W)$ are canonically isomorphic in $\D\Z$. \bb 

(R+) : For all $V,W\in\D^+\mc{G}_2$ the complexes $\RHom_{\mc{G}_1}(V,W$) and 
$\RHom_{\mc{G}_2}(V,W)$ are canonically isomorphic in $\D\Z$. \bb 

\ind Let $A$ be a commutative ring, let $Y$ be a set of prime ideals 
of $A$, let $\mc{G}_1$ (resp. $\mc{G}_2$) be the category of $A$-modules 
(resp. of $A$-modules supported on $Y$). Do (R) or (R+) hold? (See 
Theorem~\ref{qc} below for a partial answer.) \bb

\ind By the proof of Weibel \cite{weibcy}, Theorem A3, (R) implies (R+) 
\footnote
{I know no cases where (R+) holds but (R) doesn't.}. 
Moreover, if (R) (resp. (R+)) holds, then $\Ri\iota$ (resp. $\Ri^+\iota$) 
is a full embedding. Indeed we have 
$\Hom_{\D\mc{G}_i}=\HH^0\RHom_{\mc{G}_i}$ (resp. 
$\Hom_{\D^+\mc{G}_i}=\HH^0\RHom_{\mc{G}_i}$) by Lipman \cite{lip}, 
I.2.4.2. \bb 

\ind Let $\Mod A$ denote the category of left $A$-modules (whenever this 
makes sense), and let $\D A$ (resp. $\D^+A$, resp. $\K A$, resp. $\K^+A$) 
be an abbreviation for $\D\Mod A$ (resp. $\D^+\Mod A$, resp. $\K\Mod A$, 
resp. $\K^+\Mod A$), where $\K$ means ``homotopy category''. 
(Even if $\mc{G}_1$ or $\mc{G}_2$ is {\bf not} Grothendieck, it may still 
happen that (R+) or (R) makes sense and holds. In such a situation 
the phrase ``(R+) (resp. (R)) holds'' shall mean 
``(R+) (resp. (R)) makes sense and holds''.) \bb 

\ind Let $\mc{A}$ be a sheaf of rings over a topological 
space $X$, let $Y$ be a locally closed subspace of $X$, let $\mc{B}$ be the 
restriction of $\mc{A}$ to $Y$, and identify, thanks to Section 3.5 
of Grothendieck \cite{g}, $\Mod\mc{B}$ to the full subcategory of 
$\mc{A}$-modules supported on $Y$. 


\begin{ttt}\label{faisc} 
The pair $(\Mod\mc{A},\Mod\mc{B})$ satisfies (R). 
\end{ttt}

{\bf Proof.} Let $r:\Mod\mc{A}\to\Mod\mc{B}$ be the restriction 
functor. \bb 

{\it Case 1.} $Y$ is closed. --- We have for $V\in\K\mc{B}$ 

$$\Big[\Hom_\mc{B}^\pt(V,?)=\Hom_\mc{A}^\pt(V,?)\circ\K\iota\Big]:
\K\mc{B}\to\K\Z.$$

Since $\iota$ is right adjoint to the exact functor $r$, it preserves 
K-injectivity in the sense of Spaltenstein \cite{spalt}. By Lipman 
\cite{lip}, Corollaries I.2.2.7 and I.2.3.2.3, we get

$$\Big[\RHom_\mc{B}(V,?)\iso\RHom_\mc{A}(V,?)\circ\Ri\iota\Big]:
\D\mc{B}\to\D\Z.$$ 

{\it Case 2.} $Y$ is open. --- We have for $V\in\K\mc{B}$ 

$$\Big[\Hom_\mc{A}^\pt(V,?)=\Hom_\mc{B}^\pt(V,?)\circ\K r\Big]:
\K\mc{A}\to\K\Z.$$

As $r$ is right adjoint to the exact functor $\iota$, it preserves 
K-injectivity, and Lipman \cite{lip}, Corollaries I.2.2.7 and 
I.2.3.2.3, yields $\Ri r\circ\Ri\iota=\Id_{\D\mc{B}}$, 

$$\Big[\RHom_\mc{A}(V,?)=\RHom_\mc{B}(V,?)\circ\Ri r\Big]:
\D\mc{A}\to\D\Z,$$

and thus

$$\Big[\RHom_\mc{A}(V,?)\circ\Ri\iota=\RHom_\mc{B}(V,?)\Big]:
\D\mc{B}\to\D\Z.\quad\square$$

\begin{ppp}
Let $X$ and $\mc{A}$ be as above, let $Y$ be a union of closed 
subspaces of $X$, and let $\Mod(\mc{A},Y)$ be the category of 
$\mc{A}$-modules supported on $Y$. Then the pair 
$(\Mod\mc{A},\Mod(\mc{A},Y))$ satisfies (R+). 
\end{ppp} 

{\bf Proof.} See Grothendieck \cite{g}, Proposition 3.1.2, and Hartshorne 
\cite{gh}, Proposition I.5.4. $\square$ \bb 

\ind Let $(X,\mc{O}_X)$ be a noetherian scheme, $\mc{A}$ a sheaf of 
rings over $X$ and $\mc{O}_X\to\mc{A}$ a morphism, assume 
$\mc{A}$ is $\mc{O}_X$-coherent, let $Y$ be a subspace of $X$, 
and denote by $\QC\mc{A}$ (resp. $\QC(\mc{A},Y)$) 
the category of $\mc{O}_X$-quasi-coherent $\mc{A}$-modules 
(resp. $\mc{O}_X$-quasi-coherent $\mc{A}$-modules supported on $Y$). 


\begin{ttt} \label{qc} 
The pair $(\QC\mc{A},\QC(\mc{A},Y))$ satisfies (R+). 
If in addition $\Ext_{\QC\mc{A}}^n=0$ for $n\gg0$, then (R) 
holds\footnote{We regard $\Ext_\mc{G}^n$ as a functor defined on 
$\mc{G}^{op}\times\mc{G}$ (and of course {\bf not} on 
$\D\mc{G}^{op}\times\D\mc{G}$).}. 
\end{ttt}

\ind Let $A$ be a left noetherian ring, let $B$ be a ring, let 
$A\to B$ be a morphism, let $\mc{G}$ be a Grothendieck 
subcategory of $\Mod B$, let $(U_j)_{j\in J}$ be a family of 
generators of $\mc{G}$ which are finitely generated over $A$, and let 
$I$ be an Artin-Rees left ideal of $A$. For each $V$ in $\Mod A$ set 

\begin{equation}\label{V_I}
 V_I:=\{v\in V\ |\ I^{n(v)}v=0 \mbox{ for some } n(v)\in\N\}.
\end{equation}

Assume that $IV$ and $V_I$ belong to $\mc{G}$ whenever $V$ does. 
Let $\mc{G}_I$ be the full subcategory of $\mc{G}$ whose objects 
satisfy $V=V_I$. \bb

{\bf Example:} $\mc{G}$ is the category of $(\go{g},K)$-modules  
defined in Section 1.1.2 of Bernstein and Lunts \cite{bl}, 
$A$ is $U\go{g}$, $B$ is $U\go{g}\rtimes\C K$, $I$ is a left ideal of $A$ 
generated by $K$-invariant central elements.


\begin{ttt} \label{bl} 
The pair $(\mc{G},\mc{G}_I)$ satisfies (R+). 
If in addition $\Ext_\mc{G}^n=0$ for $n\gg0$, then (R) holds. 
In particular if $(\mc{G},\mc{G}_I)$ is as in the above Example and if 
$K$ is reductive, then (R) is fulfilled.
\end{ttt}

\begin{lem}\label{lbl}
If $E$ is an injective object of $\mc{G}$, then so is $E_I$. 
\end{lem} 

{\bf Lemma \ref{lbl} implies Theorem \ref{bl}.} 
By Theorem 1.10.1 of Grothendieck in \cite{g}, $\mc{G}$ and 
$\mc{G}_I$ have enough injectives. We have for $V\in\K^+\mc{G}_I$ 

$$\Big[\Hom_{\mc{G}_I}^\pt(V,?)=\Hom_\mc{G}^\pt(V,?)\circ\K^+\iota\Big]:
\K^+\mc{G}_I\to\K\Z$$

and thus, by Lemma \ref{lbl} and Hartshorne \cite{gh}, Proposition I.5.4.b, 

$$\Big[\RHom_{\mc{G}_I}(V,?)\iso\RHom_\mc{G}(V,?)\circ\Ri^+\iota\Big]:
\D^+\mc{G}_I\to\D\Z.$$

This proves the first sentence of the theorem. For the second one the 
argument is the same except for the fact we use Hartshorne \cite{gh}, 
proof of Corollary I.5.3.$\gamma$.b. (By the first sentence, 
$\Ext_\mc{G}^n=0$ for $n\gg0$ implies $\Ext_{\mc{G}_I}^n=0$ for $n\gg0$.) 
$\square$ \bb

{\bf Proof of Lemma \ref{lbl}.} Let 
$W\subset V$ be objects of $\mc{G}$ and $f:W\to E_I$ a morphism. We must 
extend $f$ to $g:V\to E_I$. We can assume, by the proof of 
Grothendieck \cite{g} Section 1.10 Lemma~1, (or by Stenstr\"om 
\cite{sten}, Proposition V.2.9), that $V$ 
is finitely generated over $A$. Since $W$ is also finitely generated over 
$A$, there is an $n$ such that $I^nf(W)=0$, and thus $f(I^nW)=0$. 
Choose a $k$ such that $W\cap I^kV\subset I^nW\subset\Ker f$ and set 

$$\ov{V}:=\frac{V}{I^kV}\quad,\qquad\ov{W}:=\frac{W}{W\cap I^kV}\quad.$$

Then $f$ induces a morphism $\ov{W}\to E_I$, which, by injectivity 
of $E$, extends to a morphism $\ov{V}\to E$, that in turn induces a 
morphism $\ov{V}\to E_I$, enabling us to define $g$ as the obvious
composition $V\to\ov{V}\to E_I$. $\square$ \bb 

\ind Let $\go{g}$ be a complex semisimple Lie algebra, let 
$\go{h}\subset\go{b}$ be respectively Cartan and Borel 
subalgebras of $\go{g}$, put $\go{n}:=[\go{b},\go{b}]$, say that the 
roots of $\go{h}$ in $\go{n}$ are positive, let $\mc{W}$ be the Weyl group 
equipped with the Bruhat ordering, let $\mc{O}_0$ be the category of 
those BGG-modules which have the generalized infinitesimal character of the
trivial module. The simple objects of $\mc{O}_0$ are parametrized by 
$\mc{W}$. Say that $Y\subset \mc{W}$ is an {\bf initial segment} if 
$x\le y$ and $y\in Y$ imply $x\in Y$, and that $w\in \mc{W}$ lies in 
the {\bf support} of $V\in\mc{O}_0$ if the simple object attached to $w$ 
is a subquotient of $V$. For such an initial segment $Y$ let 
$\mc{O}_Y$ be the subcategory of $\mc{O}_0$ consisting of objects 
supported on $Y\subset\mc{W}$. 


\begin{ttt} \label{cpso} 
The pair $(\mc{O}_0,\mc{O}_Y)$ satisfies (R). 
\end{ttt}

{\bf Proof.} In view of BGG \cite{bgg1} this will follow from 
Theorem~\ref{tech}. $\square$ \bb 


\ind Let $A$ be a ring, $I$ an ideal, and $B:=A/I$ the quotient ring.  


\begin{ttt}\label{cpst}
Assume that $\Ext_A^n(B,B)$ vanishes for $n>0$, and that there is a $p$ 
such that $\Ext_A^n(B,W)=0$ for all $n>p$ and all $B$-modules $W$. Then 
the pair $(\Mod A,\Mod B)$ satisfies (R). 
\end{ttt}

{\bf Proof.} \bb 

{\it Step 1} : $\Ext_A^n(B,W)=0$ for all $B$-modules $W$ and all 
$n>0$. --- By Theorem V.9.4 in Cartan-Eilenberg \cite{ce} we have 
$\Ext_A^n(B,F)=0$ for all free $B$-modules $F$ and all $n>0$. 
Suppose by contradiction there is an $n>0$ such that 
$\Ext_A^n(B,?)$ does not vanish on all $B$-modules; let $n$ be maximum 
for this property; choose a $B$-module $V$ such that 
$\Ext_A^n(B,V)\not=0$; consider an exact sequence $W\mono F\epi V$ with 
$F$ free; and observe the contradiction 
$0\not=\Ext_A^n(B,V)\iso\Ext_A^{n+1}(B,W)=0$.  \bb 

{\it Step 2} : Putting $r:=\Hom_A(B,?)$ we have 
$\Ri r\circ\Ri\iota=\Id_{\D B}$. --- 
The functor $r$, being a right adjoint, commutes with products, and, 
having an exact left adjoint, preserves injectives. 
Let $V$ be in $\D B$ and $I$ a Cartan-Eilenberg injective resolution 
(CEIR) of $V$ in $\Mod A$. By the previous step $rI$ is a CEIR of $rV=V$ in 
$\Mod B$. Weibel \cite{weibcy}, Theorem A3, implies 

\begin{itemize}
\item[(a)] the complex $\Tot^\Pi I\in\D A$, characterized by 

$$(\Tot^\Pi I)^n=\prod_{p+q=n}I^{pq},$$ 

is a K-injective resolution (see Spaltenstein \cite{spalt}) of $V$ in 
$\Mod A$, 
\item[(b)] $\Tot^\Pi rI=r\Tot^\Pi I$ is a K-injective resolution of 
$V=rV$ in $\Mod B$. 
\end{itemize}

Statement (a) yields: (c) $r\Tot^\Pi I=\Ri rV$. Then (b) and (c) imply 
that the natural morphism $V\to\Ri rV$ is a quasi-isomorphism. \bb

{\it Step 3} : (R) holds. --- See proof of Theorem \ref{faisc}, 
Case 2. $\square$

\begin{cor}\label{cpsc}
If there is a projective resolution $P=(P_n\mono\cdots \to P_1\to P_0)$ 
of $B$ by $A$-modules satisfying $\Hom_A(P_j,V)=0$ for all $B$-modules 
$V$ and all $j>0$, then pair $(\Mod A,\Mod B)$ satisfies (R). 
\end{cor}

\ind Let $A$ be a ring, $X$ a finite set and $e_\pt=(e_x)_{x\in X}$ a 
family of idempotents of $A$ satisfying 
$\sum_{x\in X}e_x=1$ and $e_xe_y=\delta_{xy}e_x$ (Kronecker delta) 
for all $x,y\in X$. \bb 

\ind The {\bf support} of an $A$-module $V$ is the set 
$\{x\in X\ |\ e_xV\not=0\}$. Let $\le$ be a partial ordering on $X$, and 
for any initial segment $Y$ put 

$$A_Y:=A\left/\sum_{x\notin Y}Ae_xA\right.,$$

so that $\Mod A_Y$ is the full subcategory of $\Mod A$ whose 
objects are supported on $Y$. (Here and in the sequel, for any ring $B$, 
we denote by $BbB$ the ideal generated by $b\in B$.) The image 
of $e_y$ in $A_Y$ will be still denoted by $e_y$. \bb

\ind Assume that, for any pair $(Y,y)$ where $Y$ is an initial segment 
and $y$ a maximal element of $Y$, the module $M_y:=A_Ye_y$ does {\bf not} 
depend on $Y$, but only on $y$. This is equivalent to the 
requirement that $A_Ye_y$ be supported on $\{x\in X\ |\ x\le y\}.$ \bb

\ind If $(V_\gamma)_{\gamma\in\Gamma}$ a family of $A$-modules, let 
$\langle V_\gamma\rangle_{\gamma\in\Gamma}$ denote the class of those 
$A$-modules which admit a finite filtration whose associated graded 
object is isomorphic to a product of members of the family. \bb

\ind Assume that, for any $x\in X$, the module $Ae_x$ belongs to 
$\langle M_y\rangle_{y\in X}$. 

\begin{ttt} \label{tech} 
The pair $(\Mod A,\Mod A_Y)$ satisfies (R). 
\end{ttt}

\ind This statement applies to the categories satisfying 
Conditions (1) to (6) in Section 3.2 of Beilinson, Ginzburg and 
Soergel \cite{bgs},  like the categories of BGG modules $\mc{O}_\lambda$ 
and $\mc{O}^\go{q}$ defined in Section 1.1 of \cite{bgs}, or more 
generally the category ${\mc{P}}(X,\mc{W})$ of perverse sheaves considered 
in Section 3.3 of \cite{bgs}. --- Because of the projectivity of 
$M_x=Ae_x$ we have 

\begin{lem}\label{swap} For any $x,y\in X$ with $x$ maximal  
there is a nonnegative integer $n$ and an exact sequence 
$(Ae_x)^n\mono Ae_y\epi V$ such that $V\in\langle M_z\rangle_{z<x}$. 
In particular $e_xV=0$. $\square$ 
\end{lem}

{\bf Proof of Theorem~\ref{tech}.} 
Assume $Y=X\backslash \{x\}$ where $x$ is maximal. 
Put $e:=e_x$, $I:=AeA$ and $B:=A_Y=A/I$. By the previous Lemma 
there is a nonnegative integer $n$ and an exact sequence 
$(Ae)^n\mono A\epi V$ with $IV=0$. Letting $J\subset A$ be the image of 
$(Ae)^n\mono A$, we have $J=IJ\subset I\subset J$, and thus $I=J$. 
In particular $I$ is $A$-projective and we have 
$\Hom_A(I,B)\simeq (eB)^n=0$. Corollary~\ref{cpsc} applies, proving 
Theorem~\ref{tech} for the particular initial segment $Y$. 
Lemma~\ref{swap} shows that $(B,Y,(e_y)_{y\in Y})$ satisfies the 
assumptions of Theorem~\ref{tech}, and an obvious induction completes the 
proof. $\square$\bb

\ind For any complex Lie algebra $\go{g}$ let $I_\go{g}$ be the 
annihilator of the trivial module in the center of the enveloping 
algebra. Using the notation and definitions of Knapp and Vogan \cite{kv}, 
let $(\go{g},K)$ be a reductive pair, let $(\go{g}',K')$ be a 
reductive subpair attached to $\theta$-stable subalgebra, let 
$\mc{R}^S:\mc{C}(\go{g}',K')\to\mc{C}(\go{g},K)$ be the cohomological 
induction functor defined in \cite{kv}, (5.3.b), and let $\mc{G}$ 
(resp. $\mc{G}'$) be the category of $(\go{g},K)$-modules on which 
$I_\go{g}$ (resp. $I_{\go{g}'}$) acts locally nilpotently. By \cite{kv}, 
Theorem 11.225, the functor $\mc{R}^S$ maps $\mc{G}'$ to $\mc{G}$. Let 
$F:\mc{G}'\to\mc{G}$ be 
the induced functor. By \cite{kv}, Theorem 3.35.b, $F$ is exact. It would 
be interesting to know if $F$ satisfies Condition (R). \bb

\ind Thank you to Anton Deitmar, Bernhard Keller and Wolfgang Soergel for 
their interest,  and to Martin Olbrich for having pointed out some 
mistakes in a previous version. \bb%
\newpage 


\cl{\bf Proof of Theorem~\ref{qc}} \bb 

\ind Put $\mc{O}:=\mc{O}_X$ and consider the following statements:\bb

(a) Every object of $\QC(\mc{O},Y)$ is contained into an object of 
$\QC(\mc{O},Y)$ which is injective in $\QC\mc{O}$. \bb

(b) Every object of $\QC(\mc{A},Y)$ is contained into an object of 
$\QC(\mc{A},Y)$ which is injective in $\QC\mc{A}$. \bb 

\ind We claim (a) $\then$ (b) $\then$ Theorem \ref{qc}. \bb

(a) $\then$ (b) : The functor $\mc{H}om_\mc{O}(\mc{A},?)$ preserves the 
following properties: \bb 

\ind $\pt$ quasi-coherence (by EGA I \cite{ega}, Corollary 2.2.2.vi),  

\ind $\pt$ the fact of being supported on $Y$ (by Grothendieck \cite{g}, 
 Proposition 4.1.1), 

\ind $\pt$ injectivity (by having an exact left adjoint). $\square$ \bb 

(b) $\then$ Theorem \ref{qc} : See proof of Theorem \ref{bl}.  $\square$ \bb

{\bf Proof of (a).} Let $M$ be in $\QC(\mc{O},Y)$ and let us show that 
$M$ is contained into an object of $\QC(\mc{O},Y)$ which is injective in 
$\QC\mc{O}$. We may, and will, assume that $Y$ is precisely the 
support of $M$.  \bb 

{\it Case 1.} $M$ is coherent, $(X,\mc{O})$ is affine. --- Write $A$ for 
$\G \mc{O}$, where $\G$ is the global section functor. 
Use the equivalence $\QC\mc{O}\iso\Mod A$ set up by $\G$ to work in 
the latter category. Then $M$ ``is'' a finitely generated $A$-module, 
and $Y$ is closed by Proposition II.4.4.17 in Bourbaki \cite{bac}. 
Let $I\subset A$ be the ideal of those $f$ in $A$ which vanish on $Y$, and 
$\Mod(A,Y)$ the full subcategory of $\Mod A$ whose objects are the 
$A$-modules $V$ satisfying $V=V_I$ in the sense of Notation (\ref{V_I}). 
Corollary~2 to Proposition II.4.4.17 in Bourbaki \cite{bac} 
implies that $\G$ induces a subequivalence $\QC(\mc{O},Y)\iso\Mod(B,Y)$. 
The claim now follows from Theorem~\ref{bl}. \bb

{\it Case 2.} $M$ is coherent. --- Argue as in the proof of Corollary 
III.3.6 in Hartshorne \cite{hartbk}, using Proposition 6.7.1 of EGA I 
\cite{ega}. \bb

{\it Case 3.} General case. --- By Gabriel \cite{gab} 
Corollary~1 \S II.4 (p. 358), 
Theorem~2 \S II.6 (p. 362), and 
Theorem~1 \S VI.2 (p. 443) 
we know that every object of $\QC\mc{O}$ has 
an injective hull and that any colimit of injective objects of 
$\QC\mc{O}$ is injective.
The expression $M\prec M'$, shall mean ``$M'$ is an 
injective hull of $M$ and $M\subset M'$''.  
Let $M'$ be such a hull and $Z$ the set of pairs $(N,N')$ with 

$$N\subset M,\quad N'\subset M',\quad N\prec N',\quad\Supp(N')=\Supp(N).$$

Then $Z$, equipped with its natural ordering, is inductive. 
Let $(N,N')$ is a maximal element of $Z$ and 
suppose by contradiction $N\not=M$. 
By Corollary 6.9.9 of EGA I \cite{ega} there is a $P$ such that 
$N\subset P\subset M$, $N\not=P$, and $C:=P/N$ is coherent. 
Let $\pi:P\epi C$ be the canonical projection and 
choose $P',C'$ such that $P\prec P',C\prec C'$. 
By injectivity of $N'$ there is a map $f:P\to N'$ such that 
$[N\incl P\stackrel{f}{\to} N']=[N\incl N']$ 
(obvious notation). Consider the commuting diagram

$$\xymatrix{
N'\ \ar@{^{(}->}[r]&N'\times C'&\ C'\ar@{_{(}->}[l]\\
N \ \ar@{^{(}->}[r]\ar@{^{(}->}[u]&\ P\ar@{->>}[r]^\pi\ar[u]^{f\times\pi}
                                                &C.\ar@{^{(}->}[u]}$$

We have $\Ker(f\times \pi)=\Ker(f)\cap \Ker(\pi)=\Ker(f)\cap N=0,$ 
i.e. $g:=f\times \pi$ is monic. By injectivity of $N'\times C'$ there is a 
map $P'\to N'\times C'$ such that 
$[P\incl P'\to N'\times C']=[P\stackrel{g}{\mono}N'\times C'],$ 
this map being monic by essentiality of $P\subset P'$; in particular 

$$\Supp(P')\subset \Supp(N')\cup \Supp(C').$$

A similar argument shows the existence of a monomorphism $P'\mono M'$ such 
that 
$[P\incl P'\mono M']=[P\incl M\incl M'],$ 
meaning that we can assume $P'\subset M'$. Since $(P,P')\notin Z$, this 
implies $\Supp(P)\not=\Supp(P')$, and the equalities 

$$\Supp(N')=\Supp(N)\quad (\mbox{because }\ (N,N')\in Z),$$

$$\Supp(C')=\Supp(C)\quad (\mbox{by Case~2}),$$ 

yield the contradiction

$$\Supp(P')\subset\Supp(N)\cup \Supp(C)=\Supp(P)\subset\Supp(P').\
\square$$


\cl{\bf Appendix 1}\bb

Let $k$ be a field and $\go{g}$ a Lie $k$-algebra. For $X,Y\in\D k$ put

$$\g X,Y\h:=\Hom_k^\pt(X,Y).$$

Let $C:=U\go{g}\otimes\Wedge\go{g}$ be the Koszul complex 
viewed as a differential graded coalgebra (here and in the sequel 
tensor products are taken over $k$).\bb

\ind In view of Weibel \cite{weibcy}, Theorem A3, we can define 
$\RHom_\go{g}$ by setting 

$$\RHom_\go{g}(X,Y):=\g\g C,X\h,\g C,Y\h\h^\go{g}.$$

(As usual the superscript $\go{g}$ means ``$\go{g}$-invariants''.) 
Recall that the Chevalley-Eilenberg complex, used to compute the 
cohomology of $\go{g}$ with values in $\g X,Y\h$, is defined by 
$\CE(X,Y):=\g C,\g X,Y\h\h^\go{g},$ and that there is a canonical 
isomorphism $F:\CE\iso\RHom_\go{g}$. Let 

$$\ext_{X,Y,Z}:\CE(Y,Z)\otimes\CE(X,Y)\to\CE(X,Z)$$

be the exterior product and 

$$\comp_{X,Y,Z}:\RHom_\go{g}(Y,Z)\otimes\RHom_\go{g}(X,Y)\to
\RHom_\go{g}(X,Z)$$

the composition. Then the expected formula 

$$\comp_{X,Y,Z}\ \circ\ (F_{Y,Z}\otimes F_{X,Y})
=F_{X,Z}\ \circ\ \ext_{X,Y,Z}$$

is easy to check. \bb

\cl{\bf Appendix 2}\bb

The following fact is used in various places (see for instance the 
proofs of Theorem I.3.3 in Cartan-Eilenberg \cite{ce}, Theorem 1.10.1 in 
Grothendieck \cite{g} and Lemma 4.3 in Spaltenstein \cite{spalt}). We use 
the notation and definitions of Jech \cite{jech}. 

\begin{lem}
Let $P$ be a poset, $\alpha$ a cardinal $\ge|P|$, and $\beta$ the least 
cardinal $>\alpha$. Then every poset morphism $f:\beta\to P$ is stationary. 
\end{lem}

{\bf Proof.} We can assume $P$ is infinite and $f$ is epic. The morphism 
$g:P\to\beta$ defined by $gp:=\min f^{-1}p$ satisfies $fg=\Id_P$. Put 
$\sigma:=\sup gP$. For all $p\in P$ we have $|gp|\le gp<\beta$, implying 
$|gp|\le\alpha$ for all $p$, and $\sigma\le\beta$. 
Statement (2.4) and Theorem 8 in Jech \cite{jech} entail respectively 
$\sigma=\bigcup_{p\in P}gp$ and $|P|\,\alpha=\alpha$, from which we 
conclude $|\sigma|\le\alpha$; this forces $\sigma<\beta$, that is 
$\sigma\in\beta$. For any $\gamma\in\beta$, $\gamma>\sigma$ 
we have $f\gamma=fgf\gamma\le f\sigma\le f\gamma$. $\square$


\hfill March 7, 2004\bb

\footnotesize 
\cl{Pierre-Yves Gaillard, D\'epartement de Math\'ematiques, Universit\'e 
Nancy~1, France}

\cl{\href{http://www.iecn.u-nancy.fr/~gaillard/}%
{http://www.iecn.u-nancy.fr/\~{}gaillard/}}

\end{document}